\magnification=\magstep1
\pageno=0
\overfullrule=0pt
{1.5}

\font\eightrm=cmr8


\def\real{\mathop{{\rm I}\kern-.2em\hbox{\rm R}}\nolimits}
\def\footnoterule{\kern-3pt \hrule width 2truein \kern 2.6pt}
\def\dZ{\mathop{{\rm Z}\kern-.4em\hbox{\rm Z}}\nolimits}
\def\dC{\mathop{{\rm I}\kern-.5em\hbox{\rm C}}\nolimits}
\def\SUBSECTION#1\par{\vskip0pt plus.2\vsize\penalty-75
        \vskip0pt plus-.2\vsize\bigskip\bigskip
        \leftline{\bf #1}\nobreak\smallskip\noindent}
\def\leaderfill{\leaders\hbox to 1em{\hss.\hss.}\hfill}

     \def\hphi{{\hat\phi}}
   
     \def\bk{{\bf k}}
      
\def\tr {{\tilde r}} \def\hphi{{\hat\phi}}  
\def\brho{\mathop{\rho\kern-.45em\hbox{$\rho$}}\nolimits}    
     \def\bZ{{\bf Z}}
    
 \def\b0{{\bf 0}}    
       
    \def\bR{{\bf R}}  \def\bK{{\bf K}}
\def\bell{\mathop{\ell\kern-.47em\hbox{$\ell$}}\nolimits} 
\def\bep{\mathop{\epsilon\kern-.45em\hbox{$\epsilon$}}\nolimits}
\def\btheta{\mathop{\theta\kern-.45em\hbox{$\theta$}}\nolimits}
\def\bdelta{\mathop{\delta\kern-.45em\hbox{$\delta$}}\nolimits}
\def\blambda{\mathop{\lambda\kern-.55em\hbox{$\lambda$}}\nolimits}
\def\bnu{\mathop{\nu\kern-.50em\hbox{$\nu$}}\nolimits}
\def\bgamma{\mathop{\gamma\kern-.55em\hbox{$\gamma$}}\nolimits}
\def\bomega{\mathop{\omega\kern-.62em\hbox{$\omega$}}\nolimits}
\def\bxi{\mathop{\omega\kern-.50em\hbox{$\xi$}}\nolimits}
\def\bmu{\mathop{\mu\kern-.55em\hbox{$\mu$}}\nolimits}
\def\btau{\mathop{\tau\kern-.50em\hbox{$\tau$}}\nolimits}
\def\balpha{\mathop{\alpha\kern-.54em\hbox{$\alpha$}}\nolimits}
\def\bbeta{\mathop{\beta\kern-.55em\hbox{$\beta$}}\nolimits}
\def\sqr#1#2{{\vcenter{\vbox{\hrule height.#2pt
               \hbox{\vrule width.#2pt height#1pt \kern#1pt
                \vrule width.#2pt}
              \hrule height.#2pt}}}}
 
\baselineskip=18pt

\centerline{CHARACTERIZATIONS OF ORTHONORMAL SCALE FUNCTIONS:}
\centerline{A PROBABILISTIC APPROACH}
\centerline{by V.\ Dobri{\'c}, R.\ Gundy and P.\ Hitczenko}
\bigskip
\centerline{Lehigh University, Rutgers University and North Carolina State 
University}
\vskip 4in

\bigskip
\bigskip


\vfil\eject

\noindent
{\bf Abstract.}

The construction of a multiresolution analysis starts with the
specification of a scale function. The Fourier transform of this
function is defined by an infinite product. The convergence of this
product is usually discussed in the context of $L^2(\bR)$. Here, we
treat the convergence problem by viewing the partial products as
probabilities, converging weakly to a probability defined on an 
appropriate sequence space.  We obtain a sufficient
condition for this convergence, 
which is also necessary in the case where the scale
function is continuous. These results extend and clarify
those of A.\ Cohen, and Hern{\'a}ndez, Wang, and Weiss.  
The method also applies to more general dilation schemes that
commute with translations by $\bZ^d$.

\noindent
{\bf Introduction.}

We will say that a function $\phi(x)$,
$x\in\bR$ is a {\it scaling} (or {\it scale}) function
if $\phi(x) \in L^2(\bR)$, and

(a$'$)  the function $\hat\phi(2\xi) = m(\xi)\hat\phi(\xi)$
with $m(\xi)$ a $2\pi$-periodic function in $L^2(\bR)$;

(b$'$)  $\sum_{k\in\bZ}|\hat\phi(\xi + 2\pi k)|^2 = 1$ a.e.,
$\xi\in\bR$;

(c$'$)  $\lim_{j\to\infty} |\hat\phi(\xi/2^j)| = 1$ a.e.

\noindent
This characterization of scaling functions is given in Chapter 7 of [5].
The conditions (a$'$) and (b$'$) mean that

(a)  the function ${1\over 2} \phi\big({x\over 2}\big) = \sum m_j \phi(x-j)$
with $\sum m_j^2 < \infty$;

(b) the translates of $\phi(x-k)$, $k\in\bZ$, form an orthonormal sequence
in $L^2(\bR, dx/2\pi)$.

The condition (c$'$) is independent of (a$'$) and (b$'$). However, the
``garden variety'' scale functions are integrable, with integral one,
so that $\hat\phi(\xi)$ is continuous, and $\hat\phi(0) = 1$.
In these cases, (c$'$) is satisfied. Therefore, let us assume,
for purposes of this introduction, that $\phi$ is integrable with
integral one.

How do we construct (or recognize) such functions $\phi$?
Certain features are easily discernable.   Since $\hat\phi(\xi)$
is continuous, and $\hat\phi(0) = 1$, the 
two-scale equation (a$'$) tells us that
$\hphi(\xi) = \prod_{j=1}^\infty m(\xi/2^j)$,
so that the properties of $\phi(x)$ (or $\hphi(\xi)$) are determined by
$m(\xi)$. 
If we divide the sum in (b$'$) into two parts, according to
the parity of the indices $k$, and use the variable $2\xi$, we find
$$\eqalign{
&\sum_{k\in\bZ}|\hphi(2\xi + 2\pi k)|^2\cr
&\quad = \sum_{k\in 2\bZ}|\hphi(2\xi + 2\pi k)|^2
         + \sum_{k\in 2\bZ + 1}|\hphi(2\xi + 2\pi k)|^2\cr
&\quad = |m(\xi)|^2\sum_{j\in\bZ}|\hphi(\xi + 2\pi j)|^2
         + |m(\xi + \pi)|^2\sum_{j\in\bZ}|\hphi(\xi + \pi + 2\pi j)|^2\cr
&\quad = |m(\xi)|^2 + |m(\xi + \pi)|^2\cr
&\quad = 1, \,\,{\rm a.e.}\cr}\leqno(*)$$
Therefore, in addition to $2\pi$-periodicity, the function $m(\xi)$ must
satisfy the identity expressed in the last two lines; 
since $|\hat\phi(0)|^2 = { 1}$, it also
follows that $|m(0)|^2 = { 1}$. 
However, this identity is not sufficient to insure that (b$'$) is
satisfied. 
If (b$'$) and (c$'$) together are satisfied, then
$m(\xi)$ is called a {\it low-pass filter}.  
Many authors have considered the problem of finding sufficient (and 
necessary) conditions on $m(\xi)$ so that the infinite product
$\hat\phi(\xi)$ is a scale function.
When $m(\xi)$ is a polynomial, two such sufficient conditions have been
proposed, one by Mallat [6] and the other by Daubechies (see [3], page 182).
Mallat's condition requires
$$\inf_{|\xi|\le\pi/2} |m(\xi)| > 0$$
and Daubechies:
$$m(\xi) = \big[(1+e^{i\xi})/2\big]^N {\cal L}(\xi)$$
with
$\sup_\xi |{\cal L}(\xi)| \le 2^{N-1/2}$.
The first necessary and sufficient conditions were found by 
Cohen [2], in the case where $m(\xi)$ is a polynomial;
he later extended his considerations to the case where $m(\xi)$
is $C^1(\bR)$.
The problem for more general $m(\xi)$ was considered by Hern{\'a}ndez,
Wang, and Weiss [4]. They obtained a necessary and sufficient
condition when $|m(\xi)|$ takes the values 0 and 1. In the notes of
Chapter 7 of the recent text by Hern{\'a}ndez and Weiss [5],
the authors propose the problem of finding necessary and sufficent
conditions in the case when $m(\xi)$ is not necessarily $C^1(\bR)$.
The principle purpose of this paper is to address this question. Our
results are inspired by Cohen's ideas; however, we have translated his
ideas into probabilistic terms.
This approach seems to us to be very natural for the problem at hand,
and allows us to obtain necessary and sufficient conditions in a very
general context.
In particular, the results of Cohen and Hern{\'a}ndez, Wang, and
Weiss are unified as special cases of a general theorem.
The method applies, as well, to more general dilations schemes in
$\bR^d$. These are treated in brief, in a separate section.

To motivate the probabilistic approach, let us summarize the problem
as it is usually presented. (See, for example, Daubechies ([3],  
Chapter 6.3) or Hern{\'a}ndez and Weiss ([4], Chapter 7.4).)
Given a candidate $2\pi$ periodic function  $m(\xi)$ with
$|m(\xi)|^2 + |m(\xi + \pi)|^2 = 1$, and $|m(0)|^2 = 1$, 
we form the sequence of partial products
$$\big|\hphi_N(\xi+2\pi k)\big|^2 
= \prod_{j=1}^N\bigg|m\bigg({{\xi + 2\pi k}\over 2^j}\bigg)\bigg|^2.$$
Then $\lim_{N\to\infty}|\hphi_N(\xi+2\pi k)|^2 := |\hphi(\xi+2\pi k)|^2$;
this limit is well defined a.e.\ since\break
 $|\hphi_N(\xi+2\pi k)|^2$
decreases with increasing $N$. The convergence in $L^2(\bR)$ is a 
different matter since the function $\hphi_N(\xi)$ is
$(2\pi)2^N$-periodic. Therefore, except for trivial cases, 
$|\hphi_N(\xi)|^2$ is never integrable as a function of $\xi\in\bR$.
An obvious remedy for this defect is to restrict $\hphi_N(\xi)$ to the
interval $[-2^N\pi, 2^N\pi]$. Thus, if
$$\hphi_N^*(\xi) =\cases{\hphi_N(\xi)&if $|\xi|\le 2^N\pi$\cr
                          0&otherwise,\cr}$$
then $|\hphi_N^*(\xi)|^2$ also converges to $|\hphi(\xi)|^2$ pointwise a.e.
To verify property (b$'$) in the definition of 
a scale function it turns out that it  is enough to show that
$\hphi_N^*(\xi)$ also converges in $L^2(\bR)$.
Here matters become delicate. The $L^2$ convergence is complicated by
the fact that there is no obvious domination. This is the point where
Cohen's ideas come into play. He suggested that one should modify
$\hphi_N(\xi)$ by multiplying by 
$\chi_{\lower 2pt\hbox{$_{\bK}$}}(\xi)$, 
rather than
$\chi_{[-\pi,\pi]}(\xi)$, where $\bK$ is a finite union of intervals
forming a compact set that is {\it congruent} to $[-\pi,\pi]$ in a
sense described below.
When such a $\bK$ exists, the sequence
$\hphi_N^{**}(\xi) 
= \hphi_N(\xi) \cdot \chi_{\lower 2pt\hbox{$_{\bK}$}}(\xi/2^N)$
may be shown to converge in $L^2(\bR)$. With this convergence established,
the convergence in $L^2(\bR)$ of the original sequence $\hphi_N^*(\xi)$ 
may also be proved. It was this feature of Cohen's approach that provoked
our effort to find another perspective where Cohen's condition would appear
in a more transparent fashion.  For smooth $m(\xi)$, 
Cohen's condition requires that
$$\inf_{j>0} \inf_{\xi\in {\bK}} |m(\xi/2^j)| > 0$$
where ${\bK}$ is a compact set that is a finite union of intervals, one of
which contains 0 as an interior point,  
such that ${\bK}$ is congruent to $[-\pi, \pi]$ in the following sense:

(a)  the Lebesgue measure of $\bK$ is $2\pi$;

(b)  for every $\xi\in [-\pi,\pi]$, there is a $k\in \bZ$ such that
$\xi + 2\pi k \in {\bK}$.

Notice that Cohen's condition is equivalent to a restriction on the
partial products $\hphi_N(\xi)$ for $\xi\in\bK$: Since $m(\xi)$ is
smooth in a neighborhood of the origin, the partial products converge
uniformly on any compact subset of $\bR$; therefore, the condition
may be stated as
$$\inf_{N\ge 1} \inf_{\xi\in\bK} |\hphi_N(\xi)| \ge \delta > 0,$$
for some $\delta > 0$. In fact, a more succinct way to formulate
the condition would be to omit the mention of $\bK$ altogether.  
As we shall see,  what is important is the existence of 
a lower bound $\delta$ for the infinite product. Furthermore,
it is not the topological, but the measure-theoretic character
of $\bK$ that is important: it is enough to require the lower
bound to hold almost everywhere in the following sense:

\noindent
Either

$(1)\quad |\hphi(\xi)| \ge\delta > 0$ almost everywhere in $[0,2\pi]$

\noindent
or

$(2)\quad \sup_{k\in\bZ} |\hphi(\xi + 2\pi k)| \ge\delta > 0$ 
almost everywhere in $[0,2\pi]$.

\noindent
Now let us drop the requirement that $\phi(x)$ is integrable,
and focus on properties (a$'$) and (b$'$) for a function in $L^2(\bR)$.
Given $m(\xi)$ as specified above, the conditions (1) and (2) are
sufficient for the function $\hat\phi(\xi)$, satisfying (a$'$) to
also satisfy (b$'$). When $\hphi(\xi)$ is
continuous, then conditions (1)-(2) are necessary and sufficient
for $\hphi(\xi)$ to satisfy $({\rm b}')$ {\it everywhere}, rather
than almost everywhere. In fact, we show that it is possible for
$\hphi(\xi)$ to be continuous, such that $({\rm b}')$ holds
except at two points. Furthermore, {\it any} example where
$({\rm b}')$ holds almost everywhere, but not everywhere,
is such that (1)-(2) fails for any $\delta > 0$.

The authors would like to express their gratitude to Professors G.\ Weiss and 
H.\ {\v{S}}iki{\'c}  of Washington University,  
for many valuable remarks on the 
subject of this paper. (In particular, see the section 
``The encounter$\ldots$'')

{\bf The probabilistic approach.}
In summary, we interpret the
function $|m(\xi)|^2$ as a conditional probability defined on a space of
infinite sequences. The partial products define a consistent family
of probabilities on this sequence space which converge, in the usual
(Kolmogorov) sense, to a probability.
The existence of a scale function is equivalent to the ``tightness'' of this
family of probabilities on ``finite'' sequences.

{\it The probability space.} 
Let $M(\xi) = |m(2\pi\xi)|^2$.
Notice that $M(\xi)$ is a one-periodic function that
satisfies $M(\xi) + M(\xi + 1/2) = 1$, and $M(0) = 1$.
The basic probability space $\Omega$  for our discussion is the disjoint union
 of two spaces of infinite sequences $\omega$ with coordinates $\omega_i = 0$
or 1. We will represent elements of $\Omega$ by $\{0,1\}\times\{0,1\}^{\bf
N}$; $\Omega^+$ and $\Omega^-$ will denote sequences starting with 0 and 1,
respectively. We identify integers with a subset of $\Omega$ in the following
way.  A positive integer $k$ with dyadic expansion 
$$k =
\sum_{i=1}^{\infty} \omega_i(k)2^{i-1}$$ 
is represented by the sequence
$$(0,\omega_1(k),\omega_2(k),\dots).$$ 
The integer zero is identified with the
sequence that is identically zero. A negative integer $k$ is represented by
coefficients of dyadic expansion of $-(k+1)$ preceded by 1 (thus, for example,
the sequence $(1,0,0\dots)$ represents $-1$. We denote the sequences
corresponding to nonnegative integers as $\bZ^+$, and those corresponding to
negative integers as $\bZ^-$. Fix $k\in\bZ$ and let $\bk_N = \{\omega:\omega_i
= \omega_i(k),\, 0 \le i\le N\}$ be the $N$ dimensional
$\Omega^+$-cylinder that contains $\omega(k)$. For each $\xi\in [0,1]$ 
we define
a probability $Q_\xi^N$, $0\le\xi < 1$, on the set of all such
cylinders by the following prescription. For $0 \le k\le 2^N - 1$, we set
$$Q_\xi^N(k) = \prod_{j=1}^N M\bigg({{\xi + k}\over 2^j}\bigg).$$
We then have $$ \sum_{0\le k < 2^N}\prod_{j=1}^N
M\bigg({{\xi + k}\over 2^j}\bigg) = 1,$$
where we used the basic fact that
$M(\xi) + M(\xi + 1/2) = 1$. In the language of (conditional) probability,
$$M\bigg({{\xi + k}\over 2^j}\bigg)
= Q_\xi(\omega_j(k)\Vert \omega_{j-1},\ldots \omega_1),$$
and the above sum is computed by the standard successive
conditioning procedure.

With this interpretation of $M({{\xi + k}\over 2^j})$, we see that the
product defines a probability on cylinders of  $\Omega^+$, and that
$$Q_\xi^N(\bk_N) = Q_\xi^{N+1}(\bk_N),$$
where $\bk_N$ is the $N$-dimensional cylinder corresponding to
$0\le k \le 2^N - 1$. In order to define corresponding probabilities on
$\Omega^-$ let us consider a ``reflected" filter
$$\widetilde{M}(\xi)=M(-\xi).$$ This filter may also be used to  construct a
probability on the positive integers $0\le k < 2^N$ in the same fashion, by
setting for $0\le\eta<1$ and $0\le\ell<2^N$
$$\widetilde{Q}^N_{\eta}(\ell)=\prod^N_{j=1}\widetilde{M}\big({\eta+\ell
\over2^j}\big).$$
We now define measures $P_\xi^N$ on cylinders in $\Omega$ by setting
$$P_\xi^N(k)=\cases{Q_{\xi}^{N+1}(k),& if $0\le k<2^N$;\cr
\widetilde{Q}_{1-\xi}^{N+1}\big(-(k+1)\big),& if $-2^{N}\le k<0$.\cr}$$ 
Notice that there is a double reflection, on the function, and on
the argument, and that $P_{\xi}^N$ corresponds to $N+1$ products in $Q$'s.
This specification shows that $P_\xi^N$, $N\ge 0$, is a consistent family 
($P_\xi^N(k) = P_\xi^{N+1}(k)$ for each fixed $k$), since each of the
families $Q_\xi^N$, $N\ge 1$, and $\widetilde Q_{1-\xi}^N$, $N\ge 1$
are consistent.  To see that $P_\xi^N$ defines a probability on the
integers $-2^N \le k < 2^N$, notice that 
$$\eqalign{
\sum_{-2^N\le k < 2^N} P_\xi^N(k)
&= \sum_{0\le k<2^N} Q^{N+1}(k) 
    + \sum_{-2^N\le k<0} \widetilde Q_{1-\xi}^{N+1}\big(-(k+1)\big)\cr
&= \sum_{0\le k<2^N} Q_\xi^{N+1}(k)
    + \sum_{-2^N\le k<0} Q_\xi^{N+1}(2^{N+1} + k)\cr
&= \sum_{0\le k<2^{N+1}} Q_\xi^{N+1}(k)\cr
&= 1.\cr}$$
Therefore, $P_\xi^N$, $N\ge 1$ specifies a probability on the
$\sigma$-field generated by the cylinders.

 {\bf The encounter at Washington University.} The initial version of
this paper contained an error in the formulation of the definition of
the family $P_\xi^N(\cdot)$, $N\ge 0$. We are extremely grateful to
H.\ {\v{S}}iki{\'c}, of Washington University, for showing us this error. 
The discussion with {\v{S}}iki{\'c} occurred during a visit by one of us,
to St.\ Louis, and resulted in a radical adjustment in the definition
of $P_\xi^N$. It is remarkable that the conclusions of Theorem 2 
survived this adjustment with minimal changes.

Now, we can restate the problem concerning the existence of a
scaling function in a very succinct fashion:

T{\eightrm{HEOREM}} 1. 
{\it The function $m(\xi)$ is a low-pass filter 
for a scaling function whose Fourier transform is $\hat\phi(\xi)$ 
if and only if

(b$''$) the probability $P_\xi$ is concentrated on finite sequences for
almost every $\xi$, $0\le\xi < 1$. 
We denote this by saying  $P_\xi(\bZ) = 1$ a.e.;

(c$''$) there exists a set $L^+ \subset [0,1)$ of positive measure such that
for $\xi\in L^+$,
$$\lim_{j\to\infty} |\hat\phi\big( (\xi + k)/2^j\big)|^2 = 1$$
for all $k\ge 0$, and a set of positive measure $L^- \subset [0,1)$
such that for $\xi\in L^-$,
$$\lim_{j\to\infty} |\hat\phi\big( (\xi + k)/2^j\big)|^2 = 1$$
for all $k\le -1$.}

{\bf Proof of Theorem 1.}
If $m(\xi)$ is a low-pass filter, then $\hat\phi(\xi)$ satisfies
(b$'$) which implies (b$''$). Conversely, the condition (b$''$)
is another way of stating (b$'$).

Now we must show that (c$''$) and (c$'$) are equivalent. We use the
following proposition.

P{\eightrm{ROPOSITION}} 1.
{\it Let $f(\cdot)$ be a function defined on $\bR_+$.
For $0\le\xi < 1$, and $k\ge 0$, consider the set
$$L = \bigg\{ \xi: \lim_{j\to\infty} f\big( (\xi + k)/2^j\big) = 1 \bigg\}$$
for all $k\ge 0$. This set has measure one or zero.}

{\bf Proof of Proposition 1.}
This is a special case of the Kolmogorov zero-one law. The set
$L$ is a ``tail set'' in the sense that it is invariant under all
transformations $\xi\to (\xi + k)/2^n$ for fixed $k$ and $n$,
$0\le k < 2^n$. Such invariant sets have measure zero or one.

Now, if (c$'$)holds then (c$''$) holds.
Conversely, if the (apparently) weaker condition (c$''$) holds,
the stronger condition (c$'$) holds by Proposition 1.
That is, the sets $L^\pm$ have measure one.

{\bf Remarks.}
The formulation of the second part of Theorem 1 was inspired 
by Theorem 3.16 of Papadakis, {\v S}iki{\'c}, and Weiss [7].  
They propose a 
characterization of nonnegative periodic functions $m(\xi)$ that are 
low-pass filters; this characterization assumes that the infinite 
product $\hat \phi (\xi)$ satisfies (c$'$), and they require that the 
partial products, suitably truncated, converge in $L^2({\bf R})$ to the 
limit $\hat \phi (\xi)$.  

This requirement is equivalent to our 
(b$''$) in Theorem 1. Rather than simply assume (c$'$) as they did, we chose
to state it in the present form for the following reason.
Papadakis, {\v S}iki{\' c}, and Weiss  exhibit an example
attributed to M.\ Paluszy{\' n}ski, where (b$'$) holds (that is,
the partial products converge in $L^2(\bR)$) but (c$'$) fails.
The example is simply $\hat\phi(\xi) = \chi_{[0,1)}(\xi)$. 
Clearly, $P_\xi(\bZ^-) = 0$ for all $\xi$, $0 \le\xi < 1$, and 
the condition (c$'$) does not hold: 
$\lim_{j\to\infty} \hat\phi(\xi/2^j) = 0$ if $\xi < 0$.
Therefore, in search of  minimal conditions, one might suggest
that there exist sets of positive measure $L^+$ and $L^-$
such that $P_\xi(\bZ^+) > 0$ for $\xi\in L^+$ and $P_\xi(\bZ^-) > 0$
for $\xi\in L^-$.
These conditions are {\it implied} by (c$''$).
That is, they are necessary conditions, but, in fact,  fail to be sufficient.
If the suggested necessary condition is strengthened to $P_\xi(\bZ^+) > 0$
for almost every $\xi$, $0\le\xi < 1$, the condition fails to be necessary.
Consider the Shannon filter,
$$m(\xi) = \chi_{[0,1/4)}(\xi) + \chi_{[3/4, 1)}(\xi),$$
extended periodically. Here $\hat\phi(\xi) = \chi_{[-1/2, 1/2)}(\xi)$
and $P_\xi(\bZ^+) = P_\xi(0) = 1$, with $P_\xi(\bZ^-) = 0$ if
$0\le\xi < 1/2$; also $P_\xi(\bZ^-) = P_\xi(-1) = 1$ with
$P_\xi(\bZ^+) = 0$, if $1/2 \le \xi < 1$.
Therefore, the qualification ``on a set of positive measure'' is
necessary. 
With this qualification, the suggested condition is not sufficient
to imply (c$''$). We can perturb the Shannon filter so that on a set
of positive measure $E$, such that
$$\lim |\hat\phi(\xi/2^j)|^2 = 0$$
for $\xi\in E$, but $P_\xi(\bZ) = 1$ a.e.
(We omit the details of this example.)
The upshot of all of this is that we must have $P_\xi(\bZ^+) > 0$ 
on a set of positive measure, and $P_\xi(\bZ^-) > 0$ on a set of
positive measure, as well as an almost everywhere dyadic continuity
at zero. These two requirements are captured in condition (c$''$).

The condition (c$''$) has a probabilistic interpretation in terms of
the underlying Markov chains associated with the functions $M(\xi)$.
However, we introduced the probability notions as a tool, and our
interest in the details of the probability structure are secondary.
Therefore, we chose not to express condition (c$''$) in purely
probabilistic terms, as we did for (b$''$).

In order to show that $P_\xi(\bZ) = 1$, we will use Prokhorov's criterion
of tightness for a sequence of probability measures.

{\bf Definition:}  The sequence $P_\xi^N$ is said to be tight 
on $\bZ$ in $\Omega$, if 
for every $\epsilon > 0$, there is an $n(\epsilon) = n(\epsilon,\xi) > 0$, 
such that
$$\sum_{n(\epsilon)\le |{\bk}_{N}|}P_\xi^N(\bk_N) \le \epsilon
\,\,{\rm for\,\, all}\,\, N\ge 0.$$
Here $|{\bk}_{N}|$ is the index $i$ with largest absolute value
such that $\omega_i(k) = 1$.

In terms of the integers $k\in\bZ$, we may write this tightness condition as
$$\sum_{n(\epsilon)\le k < 2^N} P_\xi^N(k) 
+ \sum_{n(\epsilon) \le k\le 2^N} P_\xi^N (-k) \le\epsilon.$$
We note that
$$\lim_{N\to\infty} P_\xi^N(\bk_N) = P_\xi(\omega(k))
= |\hat\phi(\xi + k)|^2,$$
or less formally,
$$\lim_{N\to\infty} P_\xi^N(k) = P_\xi(k).$$
Finally, we write
$$P_\xi(\bZ) = 1$$
if $\hat\phi(\xi)$ satisfies (b$'$).

{\bf Criterion:}  $P_\xi(\bZ) = 1$ if and only if $P_\xi^N$ are tight.

We omit the details of this argument.  (See Billingsley [1].)

We are now in a position to state the principal result.

T{\eightrm{HEOREM}} 2.  (i)  {\it A sufficient condition for
$P_\xi(\bZ) = 1$ almost everywhere is the following (condition (C)):}

{\it Suppose that for almost every $\xi$, 
$0\le\xi\le 1$, there exists a $\delta > 0$ 
and an integer $k(\xi)$,  such that $P_\xi(k(\xi))\ge\delta$.}

(ii) {\it Let $\xi\mapsto P_\xi(k)$ be continuous for each
$k\in\bZ$. (In other words, $|\hat\phi(\theta)|$ is continuous
for all $\theta\in\bR$, so that $\hat\phi$ satisfies the condition
(c$'$) for a scale function.) 
If condition (C) is satisfied, then
$P_\xi(\bZ) = 1$ for every $\xi$, $0\le\xi\le 1$. (That is, 
there is no exceptional set.)  Conversely, if $P_\xi(\bZ) = 1$ 
for every $\xi$, $0\le\xi\le 1$, then condition (C) holds with
no exceptional set.} 

(iii) {\it There exists a function $M(\xi)$ infinitely differentiable
at $\xi = 0$ such that $P_\xi(k)$ is continuous for each $k$, and 
such that $P_\xi(\bZ) = 1$, except at two points $\xi$, $0<\xi<1$.
At these exceptional points, $P_\xi(\bZ) = 0$. In particular, 
condition (C) fails to hold for any $\delta > 0$.}

{\bf Remark 1.}  At first, 
the distinction between ``almost everywhere'' and ``everywhere'' in the
above theorem may seem somewhat fastidious. However, these 
distinctions are crucial for the following reasons. If $\hat\phi(\xi)$
is the Fourier transform of a scale function, then the equation $({\rm b}')$
holds {\it almost everywhere}. The circumstances where $({\rm b}')$ holds
{\it everywhere} are of secondary interest. In the same spirit, the
natural assumption of the theorem concerns the behavior of $P_\xi$
almost everywhere. If, however, we require $P_\xi(k)$ to be
continuous in $\xi$ for each $k$, then the sufficient condition (a.e.) gives
the conclusion $P_\xi(\bZ) = 1$ {\it everywhere}. Conversely, if
$P_\xi(\bZ) = 1$ {\it everywhere}, then the sufficient condition 
(C)  holds
{\it everywhere}. Thus, when $P_\xi(k)$ is supposed to be continuous,
the sufficient condition (C) becomes necessary, but with a blemish:
the natural necessary condition should read, ``If $P_\xi(\bZ) = 1$
almost everywhere, then   condition (C) holds almost
everywhere.''  However part (iii) states that this cannot hold in
general, even when $P_\xi(k)$ is continuous.  In particular, there are
low-pass filters of class $C^0(\bR)$ generating continuous scale
functions that do not satisfy Cohen's condition.

{\bf Remark 2.}
When $P_\xi(k)$ is continuous for each $k$, the condition (C) is 
equivalent to that given by Cohen.
Since $P_0(0) = P_1(1) = 1$ and $P_\xi(0)$ ($P_\xi(1)$) is continuous,
there are one-sided neighborhoods of zero and one such that
$P_\xi(0)\ge \delta > 0$, $0\le\xi < \alpha$, and 
$P_\xi(1) \ge \delta > 0$, $1 - \alpha \le \xi \le 1$.
In other words, $|\hat\phi(\xi)|^2 \ge\delta > 0$ for $|\xi| \le\alpha$.
Thus, the first condition for a Cohen set is satisfied. With each $\xi_0$
we can associate an interval, $\{\xi: |\xi - \xi_0| < \epsilon\}$
centered at $\xi_0$, such that $P_\xi((k(\xi_0)) \ge \delta/2$ for
every $\xi$ in the interval. Then, we find a finite subcollection
$\xi_i$, $k(\xi_i)$, $i = 0,1,\ldots,N$, such that the corresponding
union of intervals covers the unit interval. The compact set specified
by Cohen may be constructed using translations by $k_i(\xi)$, 
$i = 0,1,\ldots,N$.

Now suppose that a compact set $\bK$, with Cohen's specifications, exists.
We will show that the probabilities $P_\xi^N(\cdot)$, $N\ge 1$ 
are tight.  Choose $n(\epsilon)$ large enough so that 
$$\sum_{n(\epsilon)\le |k|} P_\xi(k) \le\epsilon.$$
Now estimate $P_\xi^N(k)$, $n(\epsilon) \le k < 2^N$ as follows:
$$P_\xi^N(k) = P_\xi^N(k + 2^{N+1} j)$$
where $j = j((\xi + k)/2^{N+1})$ is the integer such that
$j + (\xi + k)/2^{N+1} \in\bK$. If $j\ge 0$, 
$$P_\xi^N(k) \le {1\over\delta} P_\xi^N(k+2^{N+1} j)P_{(\xi + k)/2^{N+1}}(j)$$
and $n(\epsilon) \le k \le k + 2^{N+1} j$. On the other hand, if $j < 0$,
and $n(\epsilon) \le k < 2^N$, then
$$P_\xi^N(k) \le {1\over\delta} P_\xi(-(2^{N+1} |j| - k))$$
and $n(\epsilon) \le 2^{N+1} |j| - k$ if
$n(\epsilon) \le 2^N (2 |j| - 1)$.
Therefore, either
$$P_\xi^N(k) \le {1\over\delta} P_\xi(k + 2^{N+1} |j|)$$
or
$$P_\xi^N(k) \le {1\over\delta} P_\xi(-(2^{N+1} |j| - k)).$$
In the first case,
$$\sum_{n(\epsilon)\le k < 2^N}P_\xi^N(k) 
\le {1\over\delta}\sum_{n(\epsilon)\le n} P_\xi(n) \le\epsilon/\delta;$$
in the second case,
$$\sum_{n(\epsilon)\le k < 2^N} P_\xi^N(k)
\le {1\over\delta} \sum_{2^N(2|j|-1)\le |n|}P_\xi(-n).$$
Now we choose $N$ large enough so that
$2^N + 1 \ge n(\epsilon)$.
A similar argument may be made for $k < 0$, with
$n(\epsilon) \le |k| \le 2^N$. This shows that $P_\xi^N$, $N\ge N(\epsilon)$
is tight, and therefore, that the entire collection $P_\xi^N$ is tight.

{\bf Remark 3.} Hern{\'a}ndez, Wang, and Weiss [4] treated the case
where $M(\xi)$ is a $C^1(\bR)$ function, as well as the case when 
$M(\xi)$ is a 
function taking only values 0 and 1. In the latter case
$P_\xi(k)$ also takes values 0 and 1, and the condition (C) of the
theorem becomes
$$P_\xi(k) = 1\,\,\hbox{for some}\,\, k = k(\xi)$$
for almost every $\xi$, $0\le\xi\le 1$.
This condition is obviously necessary for $P_\xi(\bZ) = 1$ a.e., in this
case. Furthermore, the ``almost everywhere'' cannot be altered.

{\bf Proof of (i).}  
Suppose $P_\xi(\cdot)$ satisfies the condition (C). Let us call the
sequence of points $\{\xi': \xi' = (\xi + k)/2^N\,\, ({\rm mod}\, 1),\,
N > 0,\, k\in\bZ\}$ the {\it orbit} of $\xi \in [0,1]$.
We want all the probabilities $P_{\xi'}$ to satisfy condition (C),
where $\xi'$ belongs to the orbit of $\xi$.
The set of ``good'' points $G$, where condition (C) holds
has full measure, and the translates of $G$ by dyadically
rational points, $G_k$ also have full measure. So, we take the set
${\tilde G} = \bigcap G_k$, of full measure, of points $\xi$
that satisfy our requirement.

Now we turn to the proof of the tightness of the sequence $P_\xi^N$, for
$\xi\in {\tilde G}$. Let $\bk_N$ denote an $N$ cylinder corresponding to the
integer $k$, as specified earlier. Let $\xi' = (\xi+k)/2^{N+1}$ (mod 1) and
$k(\xi')$ be an integer  such that
$P_{\xi'}\big(k(\xi')\big) \ge\delta$. Then the $\omega$-sequence
corresponding to $k + 2^{N+1} k(\xi')$ belongs to the $N$ cylinder $\bk_N$, and
$$P_\xi\big(k+2^{N+1}  k(\xi')\big)
= P_\xi^N(\bk_N)P_{\xi'}\big(k(\xi')\big).$$
Therefore,
$$P_\xi^N(\bk_N) \le\delta^{-1}P_\xi\big(k+2^{N+1} k(\xi')\big).$$
(This estimation is simply a transcription of Cohen's calculation.)
Now, observe that the probability $P_\xi(\cdot)$ always satisfies the
condition for tightness on $\bZ$. 
That is, for $\epsilon >0$, there exists an $n = n(\epsilon, \xi)$ such that
$$\sum_{n\le|k|} P_\xi(k) \le\epsilon$$
where $|k|$ is the largest (or smallest) index in the sequence $\omega(k)$
such that $\omega_i(k) = 1$. (This is always true since $P_\xi(\bZ) \le 1$.)
Therefore
$$\eqalign{
\sum_{n\le |k|_{N}} P_\xi^N(\bk_N)
&\le \delta^{-1} \sum_{n\le |k|\le N} P_\xi\big(k+2^{N+1} k(\xi')\big)\cr
&\le \delta^{-1} \sum_{n\le |k|} P_\xi(k)\cr
&\le \delta^{-1} \cdot \epsilon.\cr}$$
This proves that the condition of part (i) is sufficient for tightness,
and so proves that $P_\xi(\bZ) = 1$ for $\xi\in G$.

{\bf Proof of (ii).}
Now we assume the condition (C) of the theorem, and that
$P_\xi(k)$ is continuous on $[0,1]$ for each $k$. We wish to show
that $P_\xi(\bZ) = 1$ for all $\xi$ in $[0,1]$. By part (i), the
equality holds almost everywhere. If it fails at some point
$\xi_1$, $0 < \xi_1 < 1$, then there must be a point $\xi_0$ 
where $P_{\xi_0}(\bZ) = 0$.  
(Consider the possibility that $P_\xi(\bZ) > 0$ for every $\xi$.
Then, for every $\xi$, there exist $k(\xi)$ such that
$P_\xi\big(k(\xi)\big) >0$. By the continuity of $\xi\to P_\xi(k)$,
the sets $\{\xi: P_\xi(k) > 0\}$ are an open cover of $[0,1]$.
Therefore, there exists a $\delta > 0$ such that $P_\xi(\bZ) \ge \delta
> 0$ for every $\xi$.  
Then the condition
of part (i) holds everywhere, and the argument given in part (i) shows
that $P_\xi(\bZ) = 1$ everywhere. Therefore, there exist points
$\xi_0$ where $P_{\xi_0}(\bZ) = 0$.)  Since $P_\xi(0)$ is continuous,
and tends to one as $\xi$ tends to zero, $P_{\xi_0}(\bZ) = 0$ if and
only if for each $k\in\bZ$, there exists $N = N(k)$ such that 
$M\big((\xi_0 + k)/2^N\big) = 0$.
This ``sudden death syndrome'' is 
inconsistent with the hypothesis of the theorem: we will prove that,
given any $\epsilon > 0$, there exists an open set of points where
$$\max_k P_\xi(k) \le\epsilon.$$
To this end, let $\xi_0$ be a point in $[0,1]$ such that $P_{\xi_0}(\bZ)
= 0$. Now consider a neighborhood of $\xi_0$,
$$N_\eta(\xi_0) = \{\xi: |\xi-\xi_0| < \eta\},$$
where $\eta$ is chosen so that
$$|M(\xi) - M(\xi')| \le\epsilon$$
for any two points $(\xi,\xi')$ such that  $|\xi - \xi'| < \eta$
(mod 1). We claim that
$$\max_k P_\xi(k) \le\epsilon$$
in the neighborhood $N_\eta(\xi_0)$. 
If $\xi = \xi_0 + \Delta$, where $|\Delta| < \eta$,
then $M\big((\xi+k)/2^N\big) \le\epsilon$ for $N = N(k)$. This implies
$$P_\xi(k) \le \epsilon$$
for every $\xi\in N_\eta(\xi_0)$.  That is, this contradicts the assumption
if we choose $\epsilon < \delta$.

Now let us prove the necessity of the condition (C). Suppose that
$P_\xi(\bZ) = 1$ for every $\xi$ in the unit interval. This implies that
there exists a finite set of  integers $\bZ_\xi$ such that
$$P_\xi(\bZ_\xi) \ge \delta(\xi) > 0.$$
By the assumption that $P_\xi(k)$ is continuous for each $k$, the fact that
$P_\xi(\bZ) = 1$ for every $\xi$, and the compactness of $[0,1]$, we can find
a  finite set of integers $\bZ_0$, independent of $\xi$, and a fixed $\delta
> 0$, such that for $\xi$, $0\le\xi\le 1$, 
$$P_\xi(\bZ_0) \ge \delta.$$
This implies that
$$\max_{k\in \bZ_0}P_\xi(k) \ge \delta/{\rm card}(\bZ_0).$$
Thus, we have shown that condition (C) holds for every $\xi$,
$0\le\xi\le 1$.

{\bf Proof of (iii).}
The statement of part (ii) would be vacuous if it were not possible to
construct a family $P_\xi$, continuous in $\xi$ for each $k$, such that
$P_\xi(\bZ) = 1$ almost everywhere, but not everywhere. 
The following is such a construction,
inspired by an example given by Cohen [2].

Let $M(\xi)$ be a continuous periodic function, with period one, such
that $M(0) = 1$ and $M(\xi)$ is infinitely differentiable in neighborhoods
of zero, and one half. The condition $M(\xi) + M(\xi + 1/2) = 1$ is
imposed, as usual. The function $M(\xi)$ is to have only three zeros in
$0 < \xi < 1$: $M(1/2) = 0$ (dictated by the usual condition), and
$M(1/6) = M(5/6) = 0$. The latter two zeros mean that $M(1/3) = M(2/3) = 1$.
At this point, we have the example due to Cohen, cited above. However, we
insist that the function $M(\xi)$ should have cusps at the points
$\xi = 1/3$ and $\xi = 2/3$, so that
$$\sum_{k=1}^\infty\big(1 - M(1/3 \pm \epsilon/2^k)\big) = \infty$$
and
$$\sum_{k=1}^\infty\big(1 - M(2/3 \pm \epsilon/2^k)\big) = \infty,$$
for any $\epsilon$, $0 < \epsilon < 1$.  (For example, we may take
$$M(\xi) \cong 1 - \big(\log(|1/3 - \xi|)\big)^{-1}$$
for $\xi$ in a neighborhood of 1/3, with a similar specification around 2/3.)

The probability $P_\xi$, constructed using this $M(\xi)$, has the following
properties:

(a)  For any integer $k\in\bZ$, $P_\xi(k)$ is continuous in $\xi$, since 
the infinite product converges uniformly in $\xi$. ($M(\xi)$ is smooth in
a neighborhood of zero.)

(b)  $P_\xi(\bZ) = 0$ at $\xi = 1/3$ and at $\xi = 2/3$. In fact,
at the point $\xi = 1/3$, $P_\xi$ is concentrated on the single
sequence $\omega$ such that $\omega_0 = 1,\, \omega_1 = 1,\, 
\omega_2 = 0,\ldots$ 
($\omega_{2i} = 0$, $\omega_{2i+1} = 1$, $i \ge 0$); 
at the point $\xi = 2/3$, $P_\xi$ is concentrated on
$ \omega = (0,0,1,0,1,\ldots)$ (i.e., $\omega_0 = 0,
\omega_{2i-1} = 0, \omega_{2i} = 1,\,\, i \ge 1$).
On the other hand, if $\xi \ne 1/3$, $\xi\ne 2/3$, the divergence
of the above sums implies that $P_\xi(\omega) = 0$ for the two
sequences described above.
(Notice that if the cusps were placed at $\xi = 0$ and $\xi = 1$, 
rather than at $\xi = 1/3$ and $\xi = 2/3$, then $P_\xi(\bZ) \equiv 0$
for all $\xi$, $0 < \xi < 1$.)

(c)  $P_\xi(\bZ) = 1$ for all other points in the unit interval.

To prove (c) we must show that the sequence $P_\xi^N$ is tight.

To ease the burden of subscript notation,  
we will denote the cylinder $\bk_N$ by $k$. 
With this convention, we must show that
$$\sum_{n(\epsilon)\le |k| \le 2^N} P_\xi^N(k) \le \epsilon$$
for some integer $n(\epsilon)$ and all $N \ge 0$. 
Now, to find the integer $n(\epsilon)$ in the definition of tightness, 
we make a finite
number of choices, starting the process by finding $m(\epsilon)$ 
such that
$$\sum_{m(\epsilon)\le |k|} P_\xi(k) \le \epsilon.$$

Now choose $\delta$ small enough so that the interval
$(1/3 - \delta, 1/3 + \delta)$ is strictly contained in $[0,1/2]$.
Let 
$$A_\delta = (1/3 - \delta, 1/3 + \delta) \cup (-1/3 - \delta, - 1/3 + \delta),
$$
and
$$A = A(\xi,\delta) = \big\{k: |k|\le 2^N,(\xi +k)/2^{N+1} \in A_\delta\big\}.
$$
Notice that if $\xi' = (\xi + k)/2^{N+1}$, $k\in\bZ$, $|k| \le 2^N$, and
$\xi' {\not\in} A_\delta$, then 
the probability $P_{\xi'}(0)$, $\xi' > 0$ (or $P_{1+\xi'}(0)$, $\xi' < 0$)
is uniformly bounded away from zero.
In the sequel, the subscripts $-1 < \xi' < 0$ are to be interpreted as
$1 + \xi'$. Thus, with this notation, we have just stated that
$$\inf_{\xi' {\not\in} A_\delta} P_{\xi'}(0) := C^{-1}(\delta) > 0.$$
Since $P_\xi(k) = P_\xi^N(k)P_{\xi'}(0)$, we may estimate as we did in
part (ii) of the proof, to obtain
$$\sum_{ {k{\not\in} A}\atop {m(\epsilon)\le |k|\le 2^N}} P_\xi^N(k) 
\le C(\delta) \sum_{ {k{\not\in} A}\atop {m(\epsilon)\le |k|}} P_\xi(k).$$
We increase $m(\epsilon)$ to $p(\epsilon)$ if necessary, so that 
$$\sum_{p(\epsilon) \le |k|} P_\xi(k) \le C^{-1}(\delta)\epsilon.$$
Therefore,
$$\sum_{ {k{\not\in} A}\atop {p(\epsilon) \le |k| \le 2^N}}P_\xi^N(k)
\le \epsilon.$$
The ``real work'' is to estimate the sum for $k\in A$,
$|k| \ge p(\epsilon)$. If $k$ satisfies these restrictions and $k > 0$,
then $\omega_0(k) = 0$ and $\omega_N(k) = 1$, $\omega_{N-1}(k) = 0$,
$\omega_{N-2}(k) = 1,\ldots$ with this alternating pattern continuing
for a least $J$ steps. The alternating pattern is dictated by the fact
that $k/2^{N+1}$ is approximately $1/3$, which has the alternating
pattern in its dyadic expansion. The fact that the approximation is
$\delta$-close ($k\in A$) means that the alternating pattern continues
for at least $J = J(\delta)$ steps, with $\omega_{N-J}(k) = 1$.

If $- 2^N \le -k < 0$, then our convention dictates that $\omega_0(-k) = 1$
and
$$\eqalign{
-k &= - \bigg(1 + \sum_{i=1}^N \omega_i(-k)2^{i-1}\bigg)\cr
   &= - \bigg(1 + \sum_{i=1}^N \omega_i(k-1) 2^{i-1}\bigg).\cr}$$
We wish to compute the probability
$$\eqalign{
P_\xi^N(-k) &= \prod_{j=1}^{N+1} M\big( (\xi - k)/2^j\big)\cr
            &= \prod_{j=1}^{N+1} M\big( (\xi + 2^{N+1} - k)/2^j\big)\cr}$$
for $-k \in A$. This restriction $-k\in A$ means that
$\omega_N(-k) = 1$, $\omega_{N-1}(-k) = 0,\ldots$
with the alternating ones and zeros continuing for at least $J = J(\delta)$
steps.

In any case, the restriction $|k| \ge p(\epsilon)$ means that
$N \ge [\log_2 p(\epsilon)] := L$, where $[x]$ is the integer part of $x$.
To prove tightness for the entire sequence $P_\xi^N$, $N\ge 1$, it suffices
to prove tightness for $P_\xi^N$, $N\ge N(\epsilon)$. Therefore, we can
restrict our attention to $N$ such that $N-J > L$.


With this pattern in mind, we can decompose $k\in A,\,\, 0 < k < 2^n$,  
into two integers:
$$k = t_\ell + b_\ell$$
where $t_\ell$ is the ``top'' of $k$
$$t_\ell = \sum_{j=\ell+1}^{N} \omega_{j}2^{j-1}$$
where the sequence $\omega_{j},\,\, j = \ell,\ldots,N$, 
is alternately 0 and 1, as specified above.
The index $\ell$ is determined by the following rule: We observe the
sequence $\omega_{N-j},\,\, j = 0,1,\ldots,\ell$ 
which alternates between 1 and 0,
starting at $\omega_{N}=1$; we stop at the index $\ell$ 
where the coefficient
$\omega_{\ell} = 1$ and the pattern is broken for coefficients smaller than
$\ell$. (Thus,  $(\omega_{\ell-1} = 0, \omega_{\ell-2} = 0)$ 
and $\omega_{\ell-1} = 1$ are the two
possibilities when $\ell > 0$.  If the pattern is not broken, then
$\ell = 0$ and $b_\ell = 0$.)  
As we remarked above $0\le\ell\le N-J$.
This means that the ``bottom'' part of $k$,
$$b_\ell = \sum_{j=1}^\ell \omega_{j}2^{j-1}$$
has arbitrary coefficients $\omega_{j}$ for $j < \ell$, and 
$\omega_{\ell} = 1$.  Also, we note that $(\xi + b_\ell)/2^{\ell+1} {\not\in}
A_\delta$ (or $b_\ell {\not\in} A$) for any $0\le\ell\le N-J$, and
$b_\ell \ge p(\epsilon)$ if $\ell > L$.

If $k\in A$, $2^N \le k < 0$ we may carry out a similar decomposition for the
positive integer $- (k+1)$. As we have noted, $k\in A$ implies that
$- (k+1)/2^{N+1}$ is approximately $1/3$. In terms of the above notation,
$$k = - (1 + b_\ell + t_\ell)$$
and
$$\eqalign{
P_\xi^N(k) &= \tilde Q_{1-\xi}^{N+1}\big( - (k+1)\big)\cr
           &= \tilde Q_{1-\xi}^{N+1} (b_\ell + t_\ell).\cr}
$$
In this way, we see that the estimation of $P_\xi^N(k)$, for
$k < 0$, may be carried out in the same way as for $k > 0$ 
by using the reflected
filter to define probabilities on nonnegative integers. 

Now suppose that $k > 0$; we may write
$$\sum_{k\in {A}} P_\xi^N(k)
   = \sum_{\ell=0}^{N-J} \sum_b P_\xi^N(b + t).$$
(Here we have omitted the subscript $\ell$, so that $b = b_\ell$,
$t = t_\ell$.)  
Write the sum on $\ell$ in two parts
$$\sum_{\ell=0}^{N-J} \sum_b P_\xi^N(b+t)
     = \sum_{\ell=0}^L \sum_b P_\xi^N(b+t) 
         + \sum_{\ell=L+1}^{N-J}\sum_b P_\xi^N(b+t).$$

To estimate the first sum, we write each term
$$P_\xi^N(b+t) = P_\xi^\ell(b) P_{(\xi + b)/2^{\ell + 1}}^{N - \ell - 1}(t')$$
where $t' = t/2^{\ell + 1}$. Notice that $t'$ is an integer, and that the
coefficients of $t'$ satisfy
$\omega_j(t') = \omega_{\ell + 1 + j}(t)$, $j = 0,1,\ldots,N-\ell-2$.
This means that $t'$ has the same pattern as $t$. 
Since the infinite sequences of alternating zeros and ones are
assigned probability zero unless $\xi = 1/3$ or $2/3$ 
(property (b)),  we have
$$P_{(\xi+b)/2^{\ell+1}}^{N-\ell - 1}(t') = o(1)$$
as $N$ tends to infinity when $\ell \le L$, uniformly in $b = b_\ell 
{\not\in} A$.  Therefore, 
$$\eqalign{
\sum_{\ell=0}^L \sum_b P_\xi^N(b+t)
&= \sum_{\ell=0}^L \sum_b P_\xi^\ell(b)P_{(\xi+b)/2^{\ell+1}}^{N-\ell-1}(t')\cr
&\le (L+1) )\cdot o(1) = o(1)\cr}$$
as $N$ tends to infinity.  Recall here that neither $L$ nor $J$ depend
on $N$.  
That is, the above sum can be made less than
$\epsilon$ if $N \ge N(\epsilon)$.  This imposes another restriction on
the $n(\epsilon)$ we are seeking, and we incorporate this into the
calculation without further mention.

Now we estimate
$$\sum_{\ell=L+1}^{N-J} \sum_b P_\xi^N(b+t)
\le \sum_{\ell=L+1}^{N-J}\sum_b P_\xi^\ell (b).$$
Recall that $b {\not\in} A$, and $p(\epsilon) \le b$ so that
$$P_\xi^\ell(b) \le C(\delta)P_\xi(b).$$
Consequently
$$
\sum_{\ell=L+1}^{N-J} \sum_b P_\xi^N(b)
\le C(\delta) \sum_{p(\epsilon)\le b} P_\xi(b)
\le \epsilon.$$
In summary, we have shown that there exists 
$n(\epsilon) = \max\big(p(\epsilon), N(\epsilon)\big)$ such that
$$\sum_{n(\epsilon)\le |k|\le 2^N} P_\xi^N(k) \le 3\epsilon$$
for all $N$.  
This is sufficient and concludes the proof of part (iii) of the theorem.
\bigskip

{\bf The multidimensional case.} The construction of scale functions
corresponding to more general dilation schemes may be accomplished in
much the same manner as described above for the case of dyadic
dilations. Cohen's criterion may be applied without essential change.
The class of dilation schemes most frequently considered are implemented
by a matrix $A$ that maps $\bZ^d$, the integer lattice, into itself. We
assume that $A$ is {\it strictly expansive} in the sense that all 
eigenvalues $\lambda_i$ are such that $|\lambda_i| > 1$. Here, a {\it scale 
function} $\phi(x)$, $x\in\bR^d$ is a function that belongs to
$L^2\big(\bR^d/(2\pi)^d\big)$ such that 

(a$'$)$\quad \hphi(A^*\xi) = m(\xi)\hphi(\xi)$

\noindent 
for $\xi\in\bR^d$, with $m(\xi)$ periodic on the $2^d$-dimensional torus
$(2\pi)^d$ and $m(0) = 1$;

(b$'$)$\quad \sum_{k\in\bZ^d} |\hphi(\xi + 2\pi k)|^2 = 1$ a.e.

\noindent
These assumptions are not enough to insure that $\phi$ corresponds to a
multiresolution analysis since $\hphi(\xi)$ is not assumed to be continuous
at zero; however, this requirement is not relevant to the present 
discussion. (See Theorem 1.7, Chapter 2 of [5].)  For a discussion of
multiresolution analyses in this generality, see Wojtaszczyk ([9], Chapter 5).
In particular, see Proposition 5.21 {\it op.\ cit.} for a statement
of Cohen's theorem.

Given a function $\hphi(\xi)$, satisfying (b$'$), we have a probability
$P_\xi(\cdot)$ defined from $\hphi(\xi)$, that is concentrated on the
lattice $\bZ^d$, for almost every $\xi$ in any set that is congruent
to $(2\pi)^d$. (Such sets are often called {\it fundamental domains} for
the action of $(2\pi)\bZ^d$ on $\bR^d$; we shall use this term also.)
The question arises: When (b$'$) holds, does $\hphi(\xi)$
correspond to a probability on a space containing $\bZ^d$ in a manner 
similar to the case when $A = 2$, acting on $\bR$? 
The ``enveloping probability space'' is certainly not canonical, and,
the construction for the case $A=2$ has an {\it ad hoc} character.
This being so, can we describe a procedure for constructing this
probability space that applies to any dilation? 
The general case presents certain technical problems associated with the
fact that we do not know of a fundamental domain that is invariant under
the action of $(A^{-1})^*$. As a consequence, we failed in our attempts
to describe a universal sequence space $\Omega$ which is independent of
$\xi$. However, if we restrict attention to the class of transformations
that are {\it similarities}, we can carry out a construction that
generalizes the case $A=2$, and looks somewhat less impromptu than
that described above. We hope that it illuminates what was done in 
that case. A {\it similarity} is a matrix $A$ such that the eigenvalues
$\lambda_i$ have constant modulus; in our case $|\lambda_i| \equiv c > 1$.
The fundamental lemma for this construction is a result due to Strichartz 
([8], Lemma 5.1).  We quote the lemma and include its proof for
completeness.

{L{\eightrm{EMMA}} 1 (Strichartz).} 
{\it Let $B$ be a strictly expansive similarity
defined on $\bR^d$, such that $B(\bZ^d) \subset \bZ^d$. Suppose that the
common value of the modulus of any eigenvalue $\lambda$ satisfies 
$|\lambda| > 1 + d^{1/2}$.
Then there exists a set of coset representatives $r_1,r_2,\ldots,r_q$ 
($q = |\lambda|^d$) for the group $\bZ^d/B(\bZ^d)$ such that every 
element $k\in\bZ^d$ has a finite expansion}
$$k = r_{i_0} + Br_{i_1} +\cdots + B^n r_{i_n}.$$

{\bf Proof.}  The choice of coset representatives is chosen as the set
$$\bZ^d \cap B \big( (-1/2, 1/2]^d\big).$$
This is possible since the unit cube, centered at the origin, is a
fundamental domain for $\bZ^d$ acting on $\bR^d$. The element
$k\in\bZ^d$ has the coset representation
$$k = r_{i_0} + Br_{i_1} +\cdots + B^{n-1}r_{i_{n-1}} + B^n \tilde r_n$$
for some $\tilde r_n \in \bZ^d$. We must show that $\tr_n\in 
B\big((-1/2,1/2]^d\big)$ for some $n\ge 0$. Since $B$ is a similarity,
$B$ maps the ball of radius $1/2$ centered at the origin, onto a ball of 
radius $|\lambda|/2$, centered at the origin, contained in
$B\big((-1/2,1/2]^d\big)$. We will prove that $\Vert\tr_n\Vert <  
|\lambda|/2$ (that is, $\tr_n$ lies in the centered ball of radius 
$|\lambda|/2$), and so is one of the coset representatives.  
Since $\Vert B\Vert = |\lambda|$ and $|r_i|\le |\lambda|d^{1/2}/2$, we have
$$\eqalign{
\Vert B^n(\tr_n)\Vert
&\le |k| + \bigg(\sum_{i=0}^{n-1}|\lambda|^i d^{1/2}\bigg)\big(|\lambda|/2\big)
\cr
&< |k| + \big[|\lambda|^n d^{1/2}\big/(|\lambda| - 1)\big]\big(|\lambda|/2\big)
.\cr}$$
Therefore, if we take $B^{-n}$ on the left-hand side, we obtain 
$$\Vert\tr_n\Vert < |k|/|\lambda|^n 
                    + \big[d^{1/2}/(|\lambda|-1)\big](|\lambda|/2),$$ 
so that $\Vert\tr_n\Vert < |\lambda|/2$
for some $n$, as we wished to show.

Armed with the above lemma, Strichartz proved the following theorem, using
the facts about tilings of $\bR^d$.

T{\eightrm{HEOREM}} 3 (Strichartz [8]). {\it 
Let $B$ be a strictly expansive similarity transformation such that
$B(\bZ^d) \subset \bZ^d$. Suppose that the (common) value of the modulus
of any eigenvalue is greater than $1+d^{1/2}$. Let $\{r_1,r_2,\ldots,r_q\}
= {\cal R}$ be the set of coset representatives specified in Lemma 1.
Then the set $T\subset \bR^d$ defined by the equation
$$B(T) = \sum_{r_i\in {\cal R}} (T + r_i)$$
tiles $\bR^d$. That is, the Lebesgue measure of
$(T + k)\cap (T + j)$ is zero if $k\ne j$ and\break $\bigcup_{k\in\bZ^d}(T+k)
= \bR^d$.}

We refer the reader to Strichartz's paper [8], and the references there,
for a proof.

Now let us consider the problem of constructing a sequence space $\Omega$,
and an embedding of $\bZ^d\mapsto\Omega$, given a strictly expansive
similarity matrix $A$ mapping $\bZ^d$ into itself, and a candidate
function $m(2\pi\xi)$, periodic with period one, for $\xi\in\bR^d$.

A necessary (but not sufficient) condition for $M(\xi) := |m(2\pi\xi)|^2$
to be associated with a scale function (that is, a function $\hphi$
satisfying (a$'$) and (b$'$)) is that
$$\sum_{i=1}^q M\big(\xi + (A^*)^{-1} r_i\big) = 1\,\, {\rm a.e.}$$
where the integers $r_i$, $i = 1,2,\ldots,q$ are coset representatives
of the group $\bZ^d/A^*(\bZ^d)$. This follows from properties (a$'$) and
(b$'$) by an argument very similar to the one given above for the case
when $A = A^* = 2$, acting on $\bZ$. Thus, for each fixed $\xi$, we have
a probability measure concentrated on $q$ points in $\bZ^d$. It is
important to note that the measure is invariant under changes of coset
representatives. That is, if $r_i$ is replaced by $\tr_i = r_i + A^*(k)$,
$i = 1,2,\ldots,q$ for some $k\in\bZ^d$, then, since $M(\xi)$ is
periodic,
$$M\big(\xi + (A^*)^{-1} r_i\big) \equiv M\big(\xi + (A^*)^{-1}\tr_i\big)$$
for $i = 1,2,\ldots,q$.

We have assumed that $A$ is a strictly expansive similarity. Although $A$ 
does not necessarily satisfy the condition of Lemma 1, that $|\lambda| > 
1 + d^{1/2}$, there is a (smallest) integer $p$ such that $A^p$ does
fulfill this condition. The subsequence of partial products 
$$P_\xi^N(k) := \prod_{j=1}^{p\cdot N} M\big( (A^*)^{-j} (\xi + k)\big),$$
where $p$ is fixed and $N = 1,2,\ldots$ defines a sequence of probabilities
on $\bZ^d$. Each of these probabilities may be considered as a probability
on a sequence space $\Omega$ whose coordinates are integers that form a
complete set of coset representatives for the group $\bZ^d/(A^*)^p(\bZ^d)$.
The parameter set containing $\xi$ is taken to be the tile $T$ generated by
$(A^*)^p$.

To be more specific, given the candidate function $M(\xi)$ we define
$\widetilde M(\xi)$ by the product
$$\widetilde M(\xi) = \prod_{j=0}^{p-1} M\big( (A^*)^j \xi\big).$$
Now set $B = (A^*)^p$ and consider $\widetilde M(\xi)$ as a candidate 
function with the dilation matrix $B$. Notice that $\widetilde M(\xi)$
is one-periodic and
$$\prod_{j=1}^\infty \widetilde M\big(B^{-j}\xi\big) 
    = \prod_{j=1}^\infty M\big( (A^*)^{-j}\xi\big).$$
We may summarize this equality be saying that $\widetilde M(\xi)$ is the square
of the modulus of a low-pass filter for $\hphi(\xi)$ corresponding to the
dilation $B^*$. The necessary condition given above for $\hphi(\xi)$ to be 
a scale function, expressed in terms of $B$ and $\widetilde M$, becomes
$$\sum_{i=1}^{q^p} \widetilde M\big(\xi + B^{-1}r_i\big) = 1\,\,{\rm a.e.}$$
where $r_i,\,\, i = 1,2,\ldots,q^p$ is any collection of coset 
representatives for the group $\bZ^d/B(\bZ^d)$. This equality follows
from its predecessor for $A^*$. Now we are in a position to specify
$\Omega$ as a sequence space with coordinates $\omega_j(k) = r_j$ where
vectors $r_j$ are the coset representatives of $\bZ^d/B(\bZ^d)$ that
appear in the expansion
$$k = r_0 + Br_1 +\cdots + B^n r_n$$
where $n = n(k)$ is the maximal exponent in the finite expansion
provided by Lemma 1.
We let $\xi$ be the generic point in the tile $T$ generated by $B$. For each
such $\xi$, the partial products
$$P_\xi^N(k) = \prod_{j=1}^N \widetilde M\big(B^{-j}(\xi + k)\big)$$
define a sequence of consistent measures on the cylinder of $\Omega$,
as described in the one dimensional case, 
and the limiting measure $P_\xi$ is defined on
the $\sigma$-field generated by the cylinders.
It is important to note that $P_\xi^N$ defines a measure concentrated
on finite sequences $\omega(k) \in \Omega$ with $\omega_j(k) = r_j$
and $\omega_{n+j}(k) \equiv 0$ for some $n$, all $j > 0$, 
defined by the expansion given in Lemma 1:
$$k = \sum_{j=0}^n B^j r_j,\quad n = n(k).$$
Furthermore, the sets
$$Z_N = \big\{k: n(k) = N\}$$
are nested ($Z_N\subset Z_{N+1}$) and $\bZ^d = \lim Z_N$. The limiting
measure $P_\xi(\bZ^d) = 1$ if and only if the sequence $P_\xi^N$ is 
{\it tight} in the sense that given $0 < \epsilon < 1$, there exists a
set $Z_{N(\epsilon)}$ such that
$$P_\xi^{N(\epsilon) + j}(Z_{N(\epsilon)}) \ge 1 - \epsilon$$
for all $j > 0$.

Cohen's condition:  There exists a compact set $\bK$ containing a
neighborhood of the origin, and congruent to $(1/2, 1/2]^d$ such that
for $\xi\in\bK$, $M\big((A^*)^{-n}\xi\big) > 0$ for all $n\ge 1$.
The following more general condition is equivalent to Cohen's
condition when $P_\xi(k)$ is continuous for each $k\in\bZ^d$:
There exists a $\delta > 0$ and $k = k(\xi) \in
\bZ^d$ such that 
$$P_\xi\big(k(\xi)\big) \ge \delta > 0\leqno({\rm Condition\,\, C})$$
for almost every $\xi\in T$.

The proof that Condition C is sufficient for tightness is similar to the
reasoning for the case $A=2$: Given $\epsilon > 0$, find
$Z_{N(\epsilon)}$ such that $P_\xi(Z_{N(\epsilon)}^c) \le\delta\epsilon$.
Then
$$\eqalign{
P_\xi^{N(\epsilon) + j}(Z_{N(\epsilon) + j} \cap Z_{N(\epsilon)}^c)
&\le \delta^{-1} P_\xi(Z_{N(\epsilon)}^c)\cr
&\le \epsilon,\cr}$$
since for $k\in Z_{N(\epsilon) + j}$ and $\ell\in\bZ^d$
$$P_\xi\big(k + B^{N(\epsilon) + j}(\ell)\big)
      = P_\xi^{N(\epsilon) + j}(k) P_{B^{-(N(\epsilon)+j)}(\xi + k)}(\ell).$$
We conclude our discussion of the multidimensional case at this point.

\vfil\eject

\centerline{\bf References}
\bigskip
\bigskip

1.  Billingsley, P. {\bf Convergence of Probability Measures},
Wiley, New York, NY, 1968.

2.  Cohen, A. Ondelettes, analyses multir{\'e}solutions, et filtres
miroir en quadrature, {\it Ann.\ Inst.\ H.\ Poincar{\'e}, Anal.\
nonlin{\'e}aire}, {7}, 439-459, 1990.

3.  Daubechies, I. {\bf Ten Lectures on Wavelets}, (CBMS-NSF regional
conference series in applied mathematics, 61) SIAM, Philadelphia, PA,
1992.

4.  Hern{\'a}ndez, E.,  Wang, X., and Weiss, G., Smoothing minimally
supported frequency wavelets: part II, 
{\it J.\ Fourier Anal.\ Appl.}, 3 (1), 23-41, 1997.

5.  Hern{\'a}ndez, E.\ and Weiss, G. {\bf A First Course on Wavelets},
CRC Press, Inc., Boca Raton, FL, 1996.

6.  Mallat, S. Multiresolution approximation and wavelets, {\it Trans.\ 
Amer.\ Math.\ Soc.}, 315, 69-88, 1989.

7. Papadakis, M., {\v S}iki{\' c}, H., and Weiss, G.,
The characterization of low-pass filters and some basic
properties of wavelets, scaling functions, and related concepts,
{\it J.\ Fourier Anal.\ and Appl.}, 5 (5), 495-521, 1999.

8.  Strichartz, R. Wavelets and self-affine tilings, {\it Constr.\ Approx.},
9, 327-346, 1993.

9.  Wojtaszczyk, P. {\bf A Mathematical Introduction to Wavelets}, 
London Math.\ Soc.\ Student Texts 37, Cambridge Univ.\ Press,
Cambridge, U.K., 1997.
\bigskip
\bigskip
\bigskip
\bigskip
\bigskip

\noindent
{\it e-mail addresses:}

\noindent
V.\ Dobric, vd00@lehigh.edu

\noindent
R.\ Gundy, gundy@rci.rutgers.edu

\noindent
P.\ Hitczenko, pawel@math.ncsu.edu

\end